\documentclass[11pt,reqno]{article}
\usepackage[all]{xy}
\usepackage{amsmath}
\usepackage{amsthm}
\usepackage{amscd}
\usepackage{amssymb}

\numberwithin{equation}{section}

\newtheorem*{Definition}{Definition}

\newtheorem{theorem}{Theorem}[section]
\newtheorem{proposition}[theorem]{Proposition}
\newtheorem{lemma}[theorem]{Lemma}
\newtheorem{corollary}[theorem]{Corollary}

\newcounter{remarks}
{\paragraph*{Remarks}\smallskip
     \begin{list}{\arabic{remarks}. }{\usecounter{remarks}%
          \setlength{\leftmargin}{0in}%
          \setlength{\rightmargin}{0in}%
          \setlength{\labelsep}{0pt}%
          \setlength{\labelwidth}{0pt}%
          \setlength{\listparindent}{0pt}%
     }
}
{
\end{list}
}

\DeclareMathOperator{\rank}{rank}

\DeclareMathOperator{\Fix}{Fix}
\DeclareMathOperator{\SL}{SL}

\DeclareMathOperator{\GL}{GL}

\DeclareMathOperator{\Isom}{Isom}

\DeclareMathOperator\congruent{\equiv}

\newcommand\R{{\mathbf R}}

\newcommand\Q{{\mathbf Q}}

\renewcommand\H{{\mathbf H}}
\newcommand\hyp{\H}
\newcommand\C{{\mathbf C}}

\newcommand\Z{{\mathbf Z}}

\renewcommand\O{{\rm O}}

\DeclareMathOperator\CAT{CAT}

\DeclareMathOperator\chr{char}
\DeclareMathOperator\minset{Minset}
\DeclareMathOperator\Minset{Minset}

\DeclareMathOperator\spn{span}
\DeclareMathOperator\FA{FA}

\begin{document}

\title{Group actions and Helly's theorem}
\author{Benson Farb\thanks{Supported in part by the NSF.  The results of
this paper were first presented in May 2000 at MSRI.}}

\maketitle

%\tableofcontents
%\newpage

\begin{abstract}
We describe a connection between the combinatorics of generators for
certain groups and the combinatorics of Helly's 1913 theorem on convex
sets.  We use this connection to prove fixed point theorems for actions
of these groups on nonpositively curved metric spaces.  These results are encoded in a 
property that we introduce called ``property $\FA_r$'', 
which reduces to Serre's property $\FA$ when $r=1$.  The method 
applies to $S$-arithmetic groups in higher $\Q$-rank, to 
simplex reflection groups (including some non-arithmetic ones), and to 
higher rank Chevalley groups over polynomial and other rings 
(for example $\SL_n(\Z[x_1,\ldots ,x_d]), n>2$).

\end{abstract}

\section{Introduction}

Serre proved \cite{Se1,Se2} that $\Gamma=\SL_3\Z$ has {\em
property} $\FA$: any $\Gamma$-action on any simplicial 
tree $T$ has a global fixed
point.  Serre's result was extended to arithmetic and $S$-arithmetic lattices in 
higher rank semisimple groups by Margulis \cite{Ma}, to groups
with property T by Alperin and Watatani \cite{Al}, and to many 
higher-rank Chevalley groups over certain commutative rings 
by Serre, Tits and Fukunaga \cite{Fu}.

These beautiful theorems unfortunately give only a rank one/higher rank
dichotomy.  While nontrivial actions of higher-rank groups on trees do
not exist at all, these groups have many nontrivial actions on
higher-dimensional simplicial complexes of nonpositive curvature (in the
$\CAT(0)$ sense).  For example, the group $\SL_n\Z[1/p]$ acts without a
global fixed point on the affine building associated to
$\SL_n(\Q_p)$; this building is 
an $(n-1)$-dimensional, nonpositively curved simplicial
complex.  Such actions were used by Quillen and Soul\'{e} (see
\cite{So}) to prove that such groups are pushouts of
``$(n-1)$-simplices of groups'' .  

Since $\SL_n\Z[1/p]$ acts (without a global fixed point) on a
nonpositively curved simplicial complex of 
dimension $n-1$ but not one of dimension $1$, 
it is natural to ask if it admits such an action in dimension 
$m$ with $1<m<n-1$. A related question: can $\SL(n,\Z[1/p])$ be written as a 
(nice) pushout of a diagram of dimension $<n-1$?  We prove below that
the answer to these questions is ``no'', but first it 
is useful to make the following:

\begin{Definition}[Property $\FA_n$]
Let $n\geq 1$.  A group $\Gamma$ is said to have {\em property $\FA_n$} if
any isometric $\Gamma$-action on any $n$-dimensional, $\CAT(0)$
cell complex $X$ has a global fixed point.  
\end{Definition}

The property $\FA_1$ corresponds with Serre's property $\FA$.  If a group
$\Gamma$ has $\FA_n$ then it has $\FA_m$ for all $m<n$.  In this paper,
all complexes $X$ are complete and have finitely many isometry types of
cells, and all actions are by isometries. We emphasize
that $X$ need not be locally finite, and the putative action is allowed
to have infinite stabilizers and to be nonfaithful.  These possibilities
are important for applications.

There is an even stronger notion which we shall find useful.

\begin{Definition}[Strong $\FA_n$]
A group $\Gamma$ is said to satisfy the {\em strong $\FA_n$ property} if 
any $\Gamma$-action on a complete $\CAT(0)$ space $X$ satisfying the following properties 
has a global fixed point.

\medskip
\noindent
{\bf \boldmath$n$-dimensionality: }$\widetilde{H}_n(Y)=0$ for all open subsets
$Y\subseteq X$.

\smallskip
\noindent
{\bf Semisimplicity: }The action of $\Gamma$ on $X$ is {\em semisimple}, i.e.\
the translation length (see below) of each $g\in \Gamma$ is realized by some $x\in X$.
\end{Definition}

\bigskip
Bridson proved (see \cite{BrH}) that any isometric action on a cell complex as above must be semisimple, so that strong $\FA_n$ implies $\FA_n$.

When $\Gamma$ is an irreducible lattice in a semisimple Lie group, the
problem of determining when a $\Gamma$-action on a singular $\CAT(0)$
space $X$ must have a global fixed point has been studied for $X$ a tree
(Margulis \cite{Ma}), $X$ a Euclidean building 
(Margulis \cite{Ma2}, Gromov-Schoen \cite{GS}), and $X$ a $\CAT(-1)$
space (Burger-Mozes \cite{BM}, Gao \cite{Ga}).  These papers have used 
either ergodic-theoretic or harmonic maps techniques.  
In this paper we introduce a new method, which reduces to Serre's \cite{Se1,Se2} in the case of trees.
One advantage to our method is that it applies to many
finitely-generated groups which are not lattices in any locally compact
topological group.  

\bigskip
\noindent
{\bf Statement of results. }Our first main results give an 
extension of the Serre-Tits-Margulis theorems.  For the definitions in the statement of the following theorem, see Subsection \ref{subsection:arithmetic} below.

\begin{theorem}[Fixed point theorem for $S$-arithmetic groups]
\label{theorem:rigidity}
Let $k$ be an algebraic number field, and let 
$G$ be an absolutely almost simple, simply-connected, connected,  
algebraic $k$-group.  
Suppose that $r=\rank_k(G)\geq 2$.  Let $S$ be a finite set of valuations of 
$k$ containing all the archimedean valuations.  
Let $\Gamma$ be an $S$-arithmetic subgroup of 
$G$.  Then $\Gamma$ has the strong property $\FA_{r-2}$
\end{theorem}

As a special case we highlight the following.

\begin{corollary}
\label{corollary:rigidity}
Let $n\geq 3$.  Then any finite index subgroup of $\SL_n\Z$ or of 
$\SL_n\Z[1/p]$ has the strong property $\FA_{n-2}$.
\end{corollary}

The dimension in Theorem \ref{theorem:rigidity} and Corollary \ref{corollary:rigidity} 
is sharp, as can be seen from the 
action of $\SL(n,\Z[1/p])$ on the affine building for $\SL(n,\Q_p)$.   For $\SL(n,\Z)$ the 
result is not sharp: an argument of Y. Shalom combined with some
elementary arguments (see also \cite{Pa}) gives that $\SL(n,\Z),
n\geq 3$ has strong $\FA_n$ for all $n$.  
We also would like to emphasize that our
method of proof uses unipotents in a strong way, so that (unlike 
many of the other methods) we can say nothing about the $k$-rank zero case.  

Our next interest is in understanding actions of Chevalley groups over finitely generated commutative rings, for example $\SL(n,\Z[x_1,\ldots ,x_d])$.  These groups are typically not lattices in any locally compact toplogical group, and so the usual methods do not seem to apply.  

For the basics of Chevalley groups, we refer the reader to \cite{St,Ste}.  Let $\Phi$ be a reduced, irreducible root system and let $R$ be a commutative ring.  Let $G(\Phi,R)$ be a Chevalley group of type $\Phi$ over $R$, and let $E(\Phi,R)$ denote its elementary subgroup.  
Our second main result extends Fukunaga's theorem \cite{Fu} from $\FA_1$ to $\FA_n$.  

\begin{theorem}[Fixed point theorem for Chevalley groups]
\label{theorem:chevalley}
Suppose $\Phi$ has rank $r\geq 2$.  
Then the group $E(\Phi,R)$ has the strong property $\FA_{r-1}$. 
\end{theorem}

Since $\SL(n,Z[x])$ surjects onto $\SL(n,\Z[1/p])$, the example above gives that 
the dimension in Theorem \ref{theorem:chevalley} is sharp.

It is a major problem in $K$-theory to determine when $G(R)=E(R)$.  This is known for many rings, for example for polynomial or Laurent series rings over a field or over $\Z$ (this is a famous theorem of Suslin).   
We thus obtain, for example, the following.

\begin{corollary}
Let $R$ be the ring of polynomials or Laurent series over $\Z$ or over a field, and let $n\geq 3$.    
Then $\SL(n,R)$ has the strong property $\FA_{n-2}$.
\end{corollary}

Note that there 
is also overlap between the above theorems: for example $SL_n(\Z[x])$ surjects
onto $\SL_n(\Z[1/p])$, so the theorem for the former group implies it for 
the latter.  However, Theorem \ref{theorem:chevalley} can be used to deduce the corresponding result for $S$-arithmetic lattices only in those $G$ which are $k$-split .

After proving Theorem \ref{theorem:chevalley}, we wondered if
$G(\Phi,R)/E(\Phi,R)$ is always solvable when $\rank(\Phi)\geq 2$, and asked R. Hazrat and N. Vavilov, 
who then answered this question affirmatively in \cite{HV}, as long as $R$ has finite Bass-Serre dimension.  Together with Theorem \ref{theorem:chevalley}, this immediately 
gives the following.  

\begin{corollary}
Let $G(\Phi,R)$ be a Chevalley group of type $\Phi$ over a finitely generated 
commutative ring $R$ with finite Bass-Serre dimension.  Suppose $\Phi$ has rank $r\geq 2$.  
Then any action of $G(R)$ on an 
$m$-dimensional space satisfying the axioms of the strong property $\FA_m$  
with $m<r$ factors through the action of the 
solvable group $G(R)/E(R)$.
\end{corollary}

Our method also applies to certain reflection groups. 

\begin{theorem}[Reflection groups]
\label{theorem:reflection}
Let $\Gamma_n$ be any discrete group generated by reflections in the sides of
a finite volume $n$-simplex in $\R^n$ or $\hyp^n$.  Then 
$\Gamma_n$ has the strong property $\FA_{n-1}$.
\end{theorem}

The dimension $n-1$ in Theorem \ref{theorem:reflection} is 
sharp.  Also, the conclusion of this theorem is never true for polyhedra
which are not simplices; groups generated by reflections in the faces of
such polyhedra admit nontrivial actions on simplicial trees.  A. Barnhill \cite{Bar} has used 
the methods introduced in this paper to classify such actions, and to begin a 
classification for actions on higher-dimensional spaces.  We will
use Theorem
\ref{theorem:reflection} to prove that the $\Gamma_n$ above, 
even non-arithmetic examples, have properties in common with arithmetic
groups.

As a final example, we would like to mention that M. Bridson has used the method introduced in this paper to prove strong fixed point theorems for automorphism groups of free groups.

\bigskip
\noindent
{\bf Integrality and finiteness. }Representations of groups with
property $\FA_n$ are quite constrained.  Replacing ``2'' by ``n'' in
6.2.22 of \cite{Se2}, and noting that the affine building attached to
$\GL_n$ of a field with discrete valuation is a nonpositively curved
$(n-1)$-complex, immediately gives the following theorem.

\begin{theorem}[Integral eigenvalues]
\label{theorem:integral}
Suppose $\Gamma$ has property $\FA_{n-1}$.  Let $\rho:\Gamma\rightarrow 
\GL(n,k)$ be 
any representation of degree $n$ over the field $k$.  Then the
eigenvalues of each of the matrices in $\rho(\Gamma)$ are integral. In
particular they are algebriac integers if $\chr(k)=0$ and are roots of
unity if $\chr(k)>0$.
\end{theorem}

The general theory of groups with this property was introduced and
studied by Bass \cite{Ba}, who called them groups of {\em integral 
$n$-representation type}. We thus have:

\begin{center}
{\it If $\Gamma$ has $\FA_{n-1}$ then $\Gamma$ is of integral
$n$-representation type.}
\end{center}

Bass proved a finiteness theorem (Proposition 5.3 of \cite{Ba}) 
for groups of integral $n$-representation type, which immediately gives 
the following.

\begin{theorem}[Finiteness of the character variety]
\label{theorem:finiteness}
Suppose $\Gamma$ has property $\FA_{n-1}$.  Then there are only finitely 
many conjugacy classes of irreducible representations of $\Gamma$ into 
$\GL_n(K)$ for any algebraically closed field $K$.
\end{theorem}

By Theorem \ref{theorem:reflection},
groups of reflections in the faces of an $(n+1)$-simplex in 
$\H^n$ satisfy the conclusions of Theorem
\ref{theorem:integral} and Theorem \ref{theorem:finiteness}.  Note that
there are non-arithmetic examples of such reflection groups acting on 
$\H^n, n=3,4,5$ (see \cite{Vi}).  
We also note that, for reflection
groups, there are examples where the finite number in the conclusion of 
Theorem \ref{theorem:finiteness} can be made arbitrary large.  

As a concrete example consider one of the reflection groups $\Gamma_n$; 
it is a lattice in $\O(n,1)$ and has property $\FA_{n-1}$.
Hence any there are only finitely many conjugacy classes of
representation $\rho:\Gamma \rightarrow O(n,m)$ for $m<n$.

\bigskip
\noindent
{\bf Acknowledgements. }It is a pleasure to thank Dave Witte Morris for useful conversations.  

\section{Group actions on $\CAT(0)$ spaces}

{\bf Standing assumptions. }Henceforth we will assume that $X$ 
satisfies the properties as in the definition of strong property $\FA_n$.  
So $X$ is a metrically complete, $n$-dimensional, $\CAT(0)$ space, 
and all actions on $X$ are isometric and semisimple.

\subsection{Some $\CAT(0)$ preliminaries}

We first recall some facts about $\CAT(0)$ spaces.  A nice reference
for background and definitions is \cite{BrH}.

For any $g\in \Isom(X)$, let $$\minset(g)=\{x\in X: d(x,gx)\leq d(y,gy)\ 
\forall y\in X\}$$ and let $\tau(g)=\inf_{x\in X}d(x,gx)$.  The isometry 
$g$ is called {\em semisimple} if the infimum in the definition of
$\tau(g)$ is realized by some $x$.  
Under our standing assumptions, every $g\in
\Isom(X)$ is semisimple, and is 
one of two possible types: {\em elliptic}, which means that the fixed
set $\Fix(g)\neq \emptyset$; or {\em hyperbolic}, which means that
$\minset(g)\neq \emptyset$ is the unique nonempty $g$-invariant set with
every point $x\in \minset(g)$ translated a distance $\tau(g)$.

For hyperbolic $g$, the set $\minset(g)$ may be realized as a union of
parallel {\em axes} $A_g$, each of which is a $g$-invariant, bi-infinite 
geodesic line in $X$ on which $g$ acts by translation by $\tau(g)$.  
For elliptic or hyperbolic $g$, the set 
$\minset(g)$ is a closed, convex subset of $X$.  

\subsection{Finite group actions}

One of the basic results on group actions on a complete 
$\CAT(0)$ space $X$ is the following lemma.

\begin{lemma}[Bruhat-Tits Fixed Point Theorem]
\label{lemma:finite}
Any finite group action on a $\CAT(0)$ space $X$ has a global fixed point.
\end{lemma}

The idea of the proof is that the orbit of any
point under the group is a bounded set in $X$, hence has a unique 
center of mass which must be fixed by the group.  For a detailed proof
see, for example, \cite{BrH}.

\subsection{General results on $\FA_n$}
\label{subsection:general}

We now record several elementary facts.  Each holds with ``$\FA_n$'' replaced by ``strong $\FA_n$''. 

\begin{enumerate}
\item If $G$ has property $\FA_n$ then $G$ has $\FA_m$ for all $m\leq
n$.  
\label{item:lessthan}

Proof: Fixing a triangulation on each $\R^n$, one can easily convert an
action on a $\CAT(0)$ $m$-complex $X$ to an action on a $CAT(0)$
$n$-complex $Y=X\times \R^{n-m}$.  A global fixed point for the second
action gives one for the first.

\item  If $G$ has $\FA_n$ then so does every quotient of $G$.

Proof: This is obvious.

\item Let $H$ be a normal subgroup of $G$. If $H$ and $G/H$ have $\FA_n$ 
then so does $G$.

Proof: Suppose $G$ acts on a $\CAT(0)$ $n$-complex $X$.  Then
$\Fix(H)$ is a nonempty convex subset of $X$, hence is a 
$\CAT(0)$ $m$-complex for $m\leq n$.  As $H$ is normal, 
$G$ acts on $\Fix(H)$ and this action factors through $G/H$, which by
hypothesis and (\ref{item:lessthan}) above has a fixed point $y\in \Fix(H)$.  
Hence $G$ fixes $y$.

\item If some finite index subgroup $H$ of $G$ has $\FA_n$, then so does
$G$.   
\label{item:index}

Proof: Suppose $G$ acts on a $\CAT(0)$ $n$-complex $X$.  
Then $Hx=x$ for some $x\in X$.  
Let $N$ be a finite index normal subgroup of $G$ which is
contained in $H$. So $N$ fixes $x$ and $G/N$ acts on $\Fix(N)$.  As
$G/N$ is finite, we are done by Lemma \ref{lemma:finite}.

\end{enumerate}

Note that the converse to (\ref{item:index}) is not true, even when $n=1$; see,
e.g. \cite{Se2}, 6.3.5.

\subsection{Actions of nilpotent groups}

In order to understand actions of more complicated groups on $\CAT(0)$
spaces, we first need to understand actions of nilpotent groups.  
Actions of abelian groups by semisimple isometries on $\CAT(0)$ spaces
are completely understood.  The key property of such group actions is
that they have an invariant flat.  Recall that a {\em flat} in $X$ is an
isometrically embedded copy of some Euclidean space $\R^m,m\geq 0,$ in $X$.

The following proposition is well known; in this form it is (except for
the last sentence, which we prove below) Theorem 7.20 of \cite{BrH}.  

\begin{proposition}[Abelian groups]
\label{proposition:abelian}
Let $\Gamma$ be a finitely generated abelian group acting by semisimple
isometries on a $\CAT(0)$ space $X$.  Then
\begin{enumerate}
\item $\Minset(\Gamma)=\cap_{\gamma\in\Gamma}\Minset(\gamma)$ is nonempty and
splits as a product $Y\times \R^n,n\geq 0$ with $Y$ nonempty and convex.
\item $\Gamma$ leaves $\Minset(\Gamma)$ invariant, preserves the product
structure, acts trivially on the first factor, and acts cocompactly by
translations on the $\R^n$ factor. Further $n\leq \rank_{\Q}\Gamma$.
\item Any $g\in\Isom(X)$ which normalizes $\Gamma$ leaves $\Minset(\Gamma)$ 
invariant and preserves the isometric splitting $Y\times \R^n$. If $g$
in fact centralizes $\Gamma$, then the induced action on the $\R^n$
factor is by translations.
\end{enumerate}
\end{proposition}

\proof
As mentioned above, we need only prove the last claim.  By item (2), the 
action of $g$ on the $\R^n$ factor 
commutes with a cocompact group of translations.  But any such isometry
$g$ of Euclidean space must itself be a translation: if $A$ is the
linear part of $g$ then we would have $A(v+b)=Av+b$ for all $v\in \R^n$
and for every $b$ in some basis, from which it follows that $A$ is the
identity. 
\endproof

Using Proposition \ref{proposition:abelian}, we can understand nilpotent
groups acting on $\CAT(0)$ spaces.  The following is likely well-known,
and is a simple variation of the
usual theorems on properly discontinuous actions; it was proven in the
case of isometry groups of trees by Serre (see \S 6.5 of \cite{Se2}).  
We include a proof
here for completeness.

\begin{proposition}[Nilpotent groups]
\label{proposition:nilpotent}
Let $N$ be a finitely-generated, torsion-free nilpotent group acting 
by semisimple isometries on a $\CAT(0)$ space $X$.  Then either 
\begin{enumerate}
\item there is an $N$-invariant flat $\R^n, n>0$ on which $N$ acts by 
translations, hence factoring through an abelian group, or 
\item $N\cdot x=x$ for some $x\in X$.
\end{enumerate}
\end{proposition}

Note that case (2) is simply case (1) with $n=0$; we separate the
cases in the statement of the proposition for pedagogical reasons.

\proof
We prove the proposition by induction on the length $d$ of the lower central
series for $N$.  When $d=1$ then 
$N$ is abelian, and this is simply a restatement of 
Proposition \ref{proposition:abelian}, the case $n=0$ giving a global fixed point. 

Now consider an arbitrary nilpotent $N$ with lower central series of
length $d\geq 2$, and assume we know the proposition for nilpotent groups with 
lower central series of length less than $d$.  

Apply Proposition \ref{proposition:abelian} to the center $Z(N)\neq N$. So 
$\Fix(\Z(N))=Y\times \R^n, n\geq 0$, where $Y$ is a nonempty 
convex subset of $X$ (hence is a complete $\CAT(0)$ space).  As $N$
centralizes $Z(N)$, it leaves $\Minset(Z(N))$ invariant, preserves the
product decomposition $Y\times \R^n$, and acts by translations on the 
$\R^n$ factor.  In other words, the $N$-action on
$\Minset(Z(N))$ is a product action $\rho=\rho_1\times \rho_2$, where
$\rho_1:N\rightarrow \Isom(Y)$ and $\rho_2:N\rightarrow \Isom(\R^n)$.  

Since $\Fix(Z(N))=Y\times \R^n$, we have $\rho_1(Z(N))={\rm Id}$, so
$\rho_1$ factors through an $N/Z(N)$ action on the $\CAT(0)$ space $Y$.
Note that $N/Z(N)$ is nilpotent of degree less than $d$, so we can 
apply the inductive hypothesis to this action.  So $\rho_1(N/Z(N))$
leaves invariant some flat $\R^m$ in $Y$ on which it acts by
translations (the case $m=0$ being that
$\rho_1(N/Z(N))$ has a global fixed point in $Y$). Hence 
$\rho(N)=\rho_1(N)\times\rho_2(N)$ leaves invariant the flat
$\R^m\times\R^n$ in $\Minset(Z(N))$, and acts on it via translations.
\endproof

Proposition \ref{proposition:nilpotent} gives the following simple
criterion for an element in a nilpotent group $N$ to have a fixed point,
and for a nilpotent group $N$ to have a global fixed point.

\begin{corollary}[Fixed points for nilpotent groups]
\label{corollary:nilpotent:global}
\label{corollary:commutator}
Let $N$ be a finitely-generated, torsion-free nilpotent group acting by 
semisimple isometries on a $\CAT(0)$ space $X$.  Then 
\begin{enumerate}
\item If $g^m\in [N,N]$ for some $m>0$, then $g$ has a fixed point.
\item If $N$ is generated by elements each of which
has a fixed point, then $N$ has a global fixed point.
\end{enumerate}
\end{corollary}

\proof
Item (1) follows from Item (1) of Proposition
\ref{proposition:nilpotent}, since the translation subgroup of $\R^n$ is abelian and so $g$ 
must be in the kernel of the $N$-action on $\R^n$.  For Item (2), note that if (1) of Proposition 
\ref{proposition:nilpotent} were to hold, then $N$ would contain some generator 
of hyperbolic type, contradicting the given.  Thus (2) of Proposition \ref{proposition:nilpotent} 
holds, and we are done.
\endproof

\section{Helly's theorem}

One of the most basic results of convexity theory is Helly's Theorem.

\begin{theorem}[Helly, 1913]
\label{theorem:helly:original}
Let $\{X_i\}$ be any finite collection of nonempty, open convex sets in
$\R^n$.  If $$X_{i_1}\cap \cdots \cap X_{i_{n+1}}\neq \emptyset$$
for every $i_1<\cdots <i_{n+1}$, then $\bigcap_i X_i\neq \emptyset$. 
\end{theorem}

We will need a generalization of this theorem to convex subsets of 
$\CAT(0)$ $n$-complexes.  The standard proof of
Helly's theorem use {\em Radon's Theorem}: any set of at least $n+2$ points in
$\R^n$ can be divided into two nonempty subsets whose convex hulls
intersect.  However, Radon's Theorem is {\em false} when $\R^n$ is
replaced by a $\CAT(0)$ $n$-complex $X$, as can be seen even 
when $n=1$ and $X$ is a tree.  

However Helly's original proof , which seems to have been almost
completely forgotten \footnote{In fact a number of the
``generalizations'' of Helly's theorem in the literature follow
immediately from Helly's original proof.}, does generalize to this
context.  To state Helly's actual result in modern terminology, recall
that a {\em homology cell} is a nonempty space $X$ whose singular
homology groups are isomorphic to those of a point, that is the reduced
homology $\widetilde{H}_q(X)=0$ for all $q\geq 0$.  Helly's original result
and proof, restated in modern terminology by
Debrunner (see Theorem 2 of
\cite{De}), is:

\begin{theorem}[Topological Helly Theorem]
\label{theorem:helly:topological}
Let $X$ be a normal topological space with the property that every
nonempty open subset $Y\subseteq X$ has $\tilde{H}_q(Y)=0$ for all $q\geq n$.  
Let $\{X_i\}$ be any finite collection of nonempty closed homology cells
in $X$.  If the intersection of any $r$ sets of this family is nonempty 
for $r\leq n+1$, and is a homology cell for $r\leq n$, then
$\bigcap_i X_i$ is a homology cell, in particular is nonempty.
\end{theorem}

Note that nonempty convex subsets of $\R^n$ are homology cells, 
as are intersections of convex sets.  Hence Theorem
\ref{theorem:helly:topological} immediately implies Theorem
\ref{theorem:helly:original}.

Although Theorem \ref{theorem:helly:topological} is stated in \cite{De}
only in the case when $X=\R^n$, the only property of $\R^n$ used in his
proof (see Lemma $B_n$ of \cite{De}) is that every nonempty open subset
$Y\subseteq X$ has $\tilde{H}_q(Y)=0$ for all $q\geq n$.  Note that
this property holds for all contractible, 
$n$-dimensional simplicial complexes, even
ones that are not locally finite.  The proof is exactly the same as the one
given in \cite{De} for $\R^n$.

Theorem \ref{theorem:helly:topological} is stated 
in \cite{De} for open sets $X_i$; however, for a normal 
topological space, e.g. a metric space, the theorem implies the version 
for closed sets $X_i$.  This is because, for a finite collection of closed
sets in a normal topological space, there exist regular neighborhoods
around each of the sets with the property that whenever two of the sets
are disjoint so are their regular neighborhoods.

While Theorem \ref{theorem:helly:topological} suffice for 
our applications, it is actually a very special case (when $X$ is a CW complex) of the following 
old theorem of Leray (see, e.g.\ \cite{Br} Theorem VII.4.4).

\begin{theorem}[Leray]
Let $\{X_\alpha\}$ be a family of non-empty subcomplexes of a CW complex $X$.  
Suppose that every non-empty, finite 
intersection of the $X_\alpha$ is acyclic.  Then
$$H_{\ast}(\cup_i X_i)=H_\ast ({\cal N}(\{X_i\}))$$
where ${\cal N}(\{X_\alpha\})$ is the nerve of the cover $\{X_\alpha\}$.   
\end{theorem}

In particular, if $\{X_i\}$ is any collection of convex subsets of an 
$n$-dimensional metric space, and if we know the nerve of the covering 
has nonvanishing $n$-th homology, then further intersections are
forced. Theorem \ref{theorem:helly:topological} is the very special case 
when the nerve is the boundary of an $n$-simplex, in which case the 
full $n$-simplex must be part of the nerve, that is the sets have a
common intersection.  Thus Leray's theorem
gives a vast generalization of Helly's original theorem.    

It should be possible to apply the techniques of this paper to a much
wider class of groups than we consider here, replacing our use of Helly's
Theorem by Leray's Theorem.  
 
\section{Combinatorics of generators of $S$-arithmetic and Chevalley groups}

\subsection{Generators for $S$-arithmetic groups}
\label{subsection:arithmetic}

We begin by recalling some definitions.  An algebraic $k$-group $G$ 
is {\em connected} if it is connected in the Zariski topology.  
Every algebraic $k$-group is an extension of a connected algebraic $k$-group by a finite group.   
$G$ is called {\em absolutely almost simple} if any proper algebraic normal subgroup of $G$ is 
trivial.  A conencted $G$ is {\em simply connected} if every central isogeny $G'\to G$ from a connected 
algebraic $k$-group $G'$ is an algebraic group isomorphism.  The special linear group and the symplectic group are examples of 
absolutely almost simple, connected, simply-connected algebraic $k$-groups.

The goal of this subsection is to prove the following.

\begin{theorem}[Combinatorics of generators]
\label{theorem:generators}
Let $k$ be an algebraic number field, and let 
$G$ be an absolutely almost simple, simply-connected, connected,  
algebraic $k$-group.  
Suppose that $r=\rank_k(G)\geq 2$.  Let $S$ be a finite set of valuations of 
$k$ containing all the archimedean valuations.  
Let $\Gamma$ be an $S$-arithmetic subgroup of 
$G$.  Then there exists a set $C=\{\Gamma_1,\ldots
,\Gamma_{r+1}\}$ of finitely generated nilpotent subgroups of 
$\Gamma$ which has the following properties:
\begin{enumerate}
\item $C$ generates a finite index subgroup of $\Gamma$.
\item Any proper subset of $C$ generates a nilpotent subgroup of
$\Gamma$.
\item There exists $m\in \Z^+$ so that for each $\Gamma_i\in C$ there is 
a nilpotent subgroup $N<\Gamma$ so that $r^m\in [N,N]$ for all $r\in \Gamma_i$.
\end{enumerate}

\end{theorem}

\proof

Fix a maximal $k$-split torus $T$ in $G$, and let $X(T)$ denote the
group of characters on $T$.  Consider the adjoint action of $T$ on 
the Lie algebra ${\mathfrak g}$ of $G$.  Let $\Phi=\Phi_k(T,G)$  be the system of 
$k$-roots with respect to $T$.  For $\alpha\in X(T)$ let
$${\mathfrak g}_\alpha=\{X\in {\mathfrak g}: Ad(t)X=\alpha(t)X \ \forall 
t\in T\}$$
As $G$ is simple $\Phi$ is irreducible.  Let $\Delta\subset \Phi$ be a
system of simple roots and let $\Phi^+$ (resp. $\Phi^-$) be the
corresponding system of positive (resp. negative) roots.  For each
$\alpha\in \Phi$, let $U_\alpha$ be the 
unique $T$-stable subgroup of $G$ having Lie
algebra the span of the root spaces $\{{\mathfrak g}_{r\alpha}:r\in\Z^+\}$.  
The group $U_\alpha$ is a unipotent 
$k$-subgroup of $G$, and is thus nilpotent.

Let ${\cal O}_S$ be the ring of $S$-integers
in $k$, and for any subgroup $H<G$ define $H({\cal O}_S):=H\cap \GL(n,{\cal O}_S)$.  For any 
ideal $\mathfrak a\neq 0$ in
${\cal O}_S$, let 
$$H({\mathfrak a})=\{x\in H({\cal O}_S): x\congruent {\rm Id} (\mbox{mod\ } 
{\mathfrak a})\}$$

We will make essential use of the following theorem of Raghunathan, also 
proved by Margulis.

\begin{theorem}[\cite{Ra}, Theorem 1.2]
\label{theorem:raghunathan}
With notation as above, let $\mathfrak a$ be any nonzero ideal of ${\cal
O}_S$.  Then the group $E({\mathfrak a})$ 
generated by $\{U_\alpha({\mathfrak a}):\alpha \in \Phi\}$
has finite index in $G({\mathfrak a})$.
\end{theorem}

Now suppose we are given an $S$-arithmetic lattice $\Gamma$ as in the
hypothesis of the theorem.  As
$\Gamma$ is $S$-arithmetic it contains $E({\mathfrak a})$ for some ${\mathfrak a}$.  By 
Theorem \ref{theorem:raghunathan}, this containment is finite index.  By 
(\ref{item:index}) on page \pageref{item:index}, it suffices to prove
the theorem for this finite index subgroup, which we will now denote by
$\Gamma$.

Let $\Gamma_\alpha=U_\alpha({\mathfrak a})$.  Then $\Gamma_\alpha$, being an $S$-arithmetic subgroup of
the $k$-group $U_\alpha$, is a finitely generated nilpotent group.  
Note that $\Gamma_\alpha$ is nilpotent since it 
is a subgroup of the nilpotent group $U_\alpha$. 

As $\Phi$ is irreducible, there is a unique longest root $\beta'$ with
respect to $\Delta$.  Let $\beta=-\beta'$. We now set
$$C=\{\Gamma_\alpha: \alpha\in \Delta\}\cup
\{\Gamma_\beta\}$$ Note that $$|C|=|\Delta|+1=r+1$$ We claim that $C$ has
the required properties.

\bigskip
\noindent
{\bf Proof of (1): }We first show that 
$C$ generates a finite index subgroup of $\Gamma=G({\mathfrak a})$.  
For roots $\alpha,\beta$, let 
$$[\alpha,\beta]=\Phi\cap \{m\alpha +n\beta :m,n\in \Z^+\}$$

We will need the following lemma, which is the lattice analog of a
well-known fact about commutators of root groups in algebraic groups.

\begin{lemma}[\cite{Ab}, Proposition 7.2.4]  
\label{lemma:abels}
For any $\alpha,\beta\in \Phi$, let $N$ denote the group generated by
$\Gamma_\alpha$ and $\Gamma_\beta$.  Then the commutator subgroup
$[N,N]$ contains a finite index subgroup of $\Gamma_\phi$ 
for each $\phi\in [\alpha,\beta]$.  
\end{lemma}

Let $<C>$ denote the smallest subgroup of $\Gamma$ containing each group 
in $C$.  
Since $$<C>\supset \Gamma_{\Delta}:=<\{\Gamma_\alpha: \alpha \in \Delta\}>$$
it follows from Lemma \ref{lemma:abels} that 
$\Gamma_\Delta$ contains a
finite index subgroup $\Gamma'_{\alpha}$ of $\Gamma_\alpha$ for each $\alpha
\in \Phi^+$.

We claim that $<C>$ also contains a finite index subgroup
$\Gamma'_{\alpha}$ of $\Gamma_\alpha$ for each $\alpha\in \Phi^-$.  
To prove this, first note that since $-\beta$
is the unique maximal root, we can write (see \cite{Hu}, 10.4 Lemma A)
$$\beta=\sum_{\alpha\in \Delta} k_\alpha \alpha \mbox{\ \ with each \ }
k_\alpha<0$$

By definition, $<C>$ contains (a finite index subgroup of)
$\Gamma_{\alpha}$ for each $\alpha \in \Delta$, so by repeatedly applying
Lemma \ref{lemma:abels} we see that $<R>$ contains a finite index
subgroup $\Gamma'_{\alpha}$ of $\Gamma_\alpha$ 
for each $\alpha \in [\beta,\Delta]\supseteq \{-\tau: \tau \in
\Delta\}$, the final containment following from the fact that each 
$k_\alpha<0$.   

By the same argument as above which showed that the group generated by
the $\Gamma'_\alpha, \alpha\in \Delta$ contains a finite index subgroup
of each $\Gamma_\alpha, \alpha\in \Phi^+$, we deduce that $<C>$ contains
a finite index subgroup of each $\Gamma_\alpha, \alpha\in \Phi^-$,
proving the claim.

Now for each
$\alpha\in \Phi$, the finite index subgroup $\Gamma'_\alpha$ of the
nilpotent group $\Gamma_\alpha$ determines an ideal ${\mathfrak a}_\alpha$
of the finitely-generated $\Z$-algebra of $S$-integers.
Let ${\mathfrak a}_\Phi$ denote the intersection in the ring of
$S$-integers of these finitely many ideals.  Then the group generated by
all of the $\Gamma'_{\alpha}, \alpha\in \Phi$ contains the group
$\Lambda({\mathfrak a}_\Phi)$ generated by all the subgroups of the $U_\alpha$ with
entries in ${\mathfrak a}_\Phi$.  Hence $<C>$ contains
$\Lambda({\mathfrak a}_\Phi)$.  By Theorem \ref{theorem:raghunathan},
$\Lambda({\mathfrak a}_\Phi)$ has finite index in $\Gamma$, and we have
shown that $C$ generates a finite index subgroup of $\Gamma$.

\bigskip
\noindent
{\bf Proof of (2): } To prove property (2), consider any proper subset
$Q$ of $C$.  As $<Q>$ is a subgroup of $<\{U_\alpha: \alpha\in Q\}>$, it
suffices to prove that this latter group is nilpotent.  As the root
groups corresponding to any collection of positive roots generate a
unipotent subgroup of $G$ (see, e.g.. \cite{Ma2},p.37), it suffices to
prove the following.

\begin{lemma}
\label{lemma:roots2}
Any subset of $|C|-1$ elements of $\Delta \cup\{\beta\}$ is a
subset of the positive roots with respect to some basis of the root
space.
\end{lemma}

\proof 
By definition this holds for the subset $\Delta$.  So fix any
$\alpha_i \in \Delta$ and consider the collection $B=(\Delta \setminus
\alpha_i)\cup \beta$, which for convenience we will write as
$\{v_1,\ldots ,v_r=\beta\}$.

Let $(\ ,\ )$ denote the inner product on the root space.  Note that $B$
is a linearly independent set of vectors as $\{v_1,\ldots ,v_{r-1}\}$ is
part of a basis and $\beta$ has a nontrivial component in the direction 
of a vector (namely $\alpha_i$) which completes that basis.
Further we have that 
$$(v_i,v_j)\leq 0 \mbox{\ for all\ }i,j$$ 
since for
$i\neq r$ we have $v_i\in \Delta$ and $v_r$ is the negative of the maximal
root and so has inner product $\leq 0$ with each element of $\Delta$. We
now claim there is a regular vector $\gamma$ such that $(\gamma,v_i)<0$
for all $i$; in other words all of the $v_i$'s lie {\em strictly} on the
same side of the hyperplane $\gamma^{\perp}$ orthogonal to $\gamma$.  

\begin{lemma}
\label{lemma:basic:linalg}
Let $A=\{v_1,\ldots ,v_r\}$ be a set of linearly independent vectors in an 
$r$-dimensional inner product space $(V,(\ ,\ ))$.  Suppose $v\in V$ is a
nonzero vector with $(v,v_i)\leq 0$ for all $i$.  Then there exists 
a vector $\gamma$ with $(\gamma,v_i)<0$ for each $i$.  Further, $\gamma$ 
can be chosen so that it 
does not lie in any subspace spanned by any proper subset of $A$.
\end{lemma}

\proof
Consider $v_j\in v^\perp$; if no such $v_i$ exists then we are done.  
Let 
$$W=\spn\{v_i: v_i\in v^\perp, i\neq j\}$$
Write $v_j=u+u'$ with $u\in \spn\{v^\perp\}\cap W$ and with $u'\in
\spn\{v^\perp\}\cap W^\perp$.  By linear independence of $A$, we have
$u'\neq 0$.  Let $\epsilon>0$, and let 
$$\gamma=v-\epsilon u'$$
Then 
$$
\begin{array}{rl}
(\gamma,v_j)&=(v-\epsilon u',u+u')\\
&=(v,u)+(v,u')-\epsilon (u',u)-\epsilon(u',u')\\
&=0+0-0-\epsilon||u'||^2\\
&<0
\end{array}$$
and for every $w\in W$
$$
\begin{array}{rl}
(\gamma,w)&=(v-\epsilon u',w)\\
&=(v,w)-\epsilon (u',w)\\
&=(v,w)
\end{array}$$

In other words, by replacing $v$ with $\gamma$ we have made the inner
product with $v_j$ strictly negative while keeping the inner product
with each $v_i\in W$ nonpositive.  Further, since there is a fixed bound 
away from zero for $(v,v_k)$ with $v_k\not\in W$, by taking $\epsilon$
sufficiently small we may guarantee that $(\gamma,v_k)<0$ for each
$v_k\not\in W$.  

Repeating the above process until $W$ is empty completes the proof of
the first claim of the lemma.  The second claim is clear.
\endproof

Lemma \ref{lemma:basic:linalg} clearly implies the claim; the last claim 
of the lemma giving that $\gamma$ is regular.

Continuing with the proof of Lemma \ref{lemma:roots2}, for any regular vector $\gamma$ let 
$$\Phi^+(\gamma)=\{\alpha\in \Phi: (\alpha,\gamma)>0\}$$
Then (see, e.g.\  Theorem 10.1 of \cite{Hu}) the set of all
indecomposable roots in $\Phi^+(\gamma)$ is a root basis of $\Phi$, so
in particular every element of $\Phi^+(\gamma)$ is a positive root with
respect to that basis.  Choosing $\gamma$ as above, we have shown that
$B\subseteq \Phi^+(\gamma)$, and we are done.
\endproof

\bigskip
\noindent
{\bf Proof of (3): }  Let $\alpha\in \Delta\cup \{\beta\}$ be given (in fact we will only use that
$\alpha\in \Phi$).  As $\Phi$ is irreducible and has rank at least $2$,
there exists $\sigma\in \Phi$ which is not proportional to $\alpha$ and
so that $(\alpha,\sigma)>0$.  Then (see, e.g.\ \cite{Hu}, Lemma 9.4)
$\alpha-\sigma\in \Phi$.  Now the group $N$ generated by
$\Gamma_{\sigma}$ and $\Gamma_{\alpha-\sigma}$ is nilpotent by the proof 
of (2), since any
two nonproportional roots lie strictly on the same 
side of the hyperplane orthogonal to some regular vector; here again we
are using that $\Phi$ is irreducible of rank at least $2$.  
Further, by Lemma \ref{lemma:abels}, the
commutator subgroup $[N,N]$ contains a finite index subgroup of
$\Gamma_{\sigma+\alpha-\sigma}=\Gamma_{\alpha}$.  In particular some
positive power of the generator $g_\alpha\in \Gamma_\alpha$ lies in
$[N,N]$.
\endproof

\subsection{Chevalley groups over commutative rings}

\bigskip
\noindent
{\bf Basics of Chevalley groups.  }We now recall the definition and basic properties of
Chevalley groups.  The reader is referred to \cite{Hur,St,Ste} for
details.  

Let $R$ be any commutative ring (with unit) which is finitely-generated
as a $\Z$-algebra.  Let $\Phi$ be a reduced, irreducible root system,  $\Delta$
a base for $\Phi$, and $\Phi^+$ the set of positive roots with respect
to $\Delta$.  Each such $\Phi$ determines a unique Lie algebra
$\mathfrak{g}(\Phi)$ over $\C$.  Fixing a faithful, 
finite-dimensional, complex representation of $\mathfrak{g}(\Phi)$
determines a {\em Chevalley-Demazure group scheme}, or {\em Chevalley
group}, $G(\Phi,R)=G_\rho(\Phi,R)$ over $R$.  We always assume that 
$G(\Phi,\C)$ is simply connected.  

For each root $\alpha\in \Phi$ there is a group isomorphism 
$$t\mapsto x_{\alpha}(t)$$
from $R$ onto a subgroup $X_\alpha$ of $G(\Phi,R)$, called the {\em root
subgroup} of $G(\Phi,R)$ corresponding to $\alpha$.  The group generated
by $\{X_\alpha: \alpha\in \Phi\}$ is called the group of {\em elementary
matrices}, and is denoted by $E(\Phi,R)$.  We also let $U(\Phi,R)$ be
the group generated by $\{X_\alpha: \alpha\in \Phi^+\}$.  Then
$U(\Phi,R)$ is nilpotent (see, e.g. \cite{Ste,St}).

There are two fudnamental formulas that hold in a Chevalley group.  The first is elementary: 
for all $s,t\in R$ and $\alpha \in \Phi$ we have

\begin{equation}
x_\alpha(s)x_\alpha(t)=x_\alpha(s+t)
\end{equation}
The second formula, which holds when $\Phi$ has rank at least $2$, 
is the {\em Chevalley Commutator Formula}: for $s,t\in R$ and for
linearly independent $\alpha,\beta \in \Phi$,
\begin{equation}
\label{eq:chevalley}
[x_\alpha(s),x_\beta(t)]=\prod_{i\alpha+j\beta \in
\Phi}x_{i\alpha+j\beta}(N_{\alpha,\beta,i,j}s^it^j)
\end{equation}
where $i,j>0$ and $N_{\alpha,\beta,i,j}$ are certain integers.

\bigskip
\noindent
{\bf A nice generating set for $E(\Phi,R)$.  }
We assume henceforth that $\Phi$ is reduced, irreducible, and has rank
at least $2$.  Let $\{r_1,\ldots ,r_n\}$ be a generating set for $R$.
As $\Phi$ is irreducible it has a unique maximal root $\beta$ (see,
e.g., Lemma 10.4.A of \cite{Hu}).  Let 
$$C=\Delta \cup \{\beta\}$$
Note that $$|C|=\rank(\Phi)+1$$
For each $\alpha \in \Delta \cup \{\beta\}$, let $\Gamma_\alpha$ denote the subgroup
of $E(\Phi,R)$ generated by $\{x_\alpha(1),x_\alpha(r_1),\ldots
,x_\alpha(r_n)\}$.

The following proposition combines several results of Fukunaga \cite{Fu}.

\begin{proposition}[Combinatorics of generators for Chevalley groups]
\label{proposition:fukunaga}
With notation as in the previous paragraph, the following holds.
\begin{enumerate}
\item $E(\Phi,R)$ is generated by $\{\Gamma_\alpha: \alpha\in C\}$.
\item For each $\alpha\in C$ and each $r\in \Gamma_\alpha$, there 
exists $p=1,2,3$ or $6$ so that $r^p\in [U,U]$ for some nilpotent 
subgroup $U\subset E(\Phi,R)$.
\item For any proper subset $C'\subset C$, the group generated by
$\{\Gamma_\alpha: \alpha\in C'\}$ is nilpotent.
\end{enumerate}
\end{proposition}

\proof
Since $\Phi$ is irreducible and $\rank(\Phi)\geq 2$, we may apply 
Proposition 2 of \cite{Fu}, which states that for each $t\in R$ and each
$\alpha\in \Phi$, there exists $p=1,2,3$ or $6$ and some base $\Delta'$
so that $x^p_\alpha(t)\in [U,U]$, where $U=U(\Phi,\Delta')$, which is
nilpotent.  This immediately gives (2).  We may also apply Proposition 4 of \cite{Fu}, which 
gives precisely (1).   The proof of item (2) of Theorem \ref{theorem:generators} above, in particular Lemma \ref{lemma:roots2}, immediately gives (3).
\endproof

\section{Finishing the proofs of the main theorems}

\subsection{A general fixed point theorem}

Our results on both $S$-arithmetic groups and Chevalley grops admit a common generalization.

\begin{theorem}
\label{theorem:general}
Let $\Gamma$ be a finitely generated group, and let $C=\{\Gamma_1,\ldots ,\Gamma_{r+1}\}$ be a collection of finitely generated nilpotent subgroups of $\Gamma$.  Suppose that:
\begin{enumerate}
\item $C$ generates a finite index subgroup of $\Gamma$.
\item Any proper subset of $C$ generates a nilpotent group.
\item There exists $m>0$ so that for any element $r$ of any $\Gamma_i$ there is a nilpotent subgroup 
$N<\Gamma$ with $r^m\in[N,N]$.
\end{enumerate}
Then $\Gamma$ has Property $\FA_{r-1}$.
\end{theorem}

\proof
Suppose that $\Gamma$ acts on an $(r-1)$-dimensional $\CAT(0)$ complex $X$.   By item \ref{item:index} in Section \ref{subsection:general},  
assumption (1) of the theorem gives that 
it suffices to prove that the subgroup generated by $C$, which we will by abuse call
$\Gamma$, fixes some point of $X$.

By assumption (3) of the theorem, we may apply part (1) of 
Corollary \ref{corollary:commutator}
to give that $\Fix(r)\neq \emptyset$ for each $r$ in each $\Gamma_i$, and then we may apply part (2) 
of that corollary to give that $\Fix(\Gamma_i)\neq\emptyset$ for each $1\leq i\leq r+1$.   
Hence we have a collection $\{\Fix(\Gamma_i)\}$ of $r+1$ 
nonempty, closed convex sets.  As each $\Fix(\Gamma_i)$ is convex it is contractible.  Since 
$C$ generates $\Gamma$ it is enough to show that the intersection of these 
$r+1$ sets is nonempty, thus giving a global fixed point for $\Gamma$.

Let $U$ be any subgroup of $\Gamma$ generated by any subset $Q\subset C$ of 
$|C|-1=r$ elements.  By assumption (2),  $U$ is
nilpotent.  By the previous paragraph, $U$ is generated by elements each of which has a fixed
point.  Corollary \ref{corollary:nilpotent:global} now implies that 
$\Fix(U)$ is a nonempty convex set.

We thus have $r+1$ contractible sets $\{\Fix(\Gamma_i)\}$ in the contractible space $X$, and the intersection of any $r$ of these sets is nonempty.  Applying the topological Helly theorem 
(Theorem \ref{theorem:helly:topological}) then gives $\bigcap_{i=1}^r \Fix(\Gamma_i)
\neq \emptyset$, and we are done.
\endproof

\subsection{Proof of Theorem \ref{theorem:rigidity}}

Item (1) of Theorem \ref{theorem:rigidity} follows immediately from 
Theorem \ref{theorem:generators} together with Theorem \ref{theorem:general}.  
Item (3) of Theorem \ref{theorem:rigidity} follows immediately from Proposition 
\ref{proposition:fukunaga} together with Theorem \ref{theorem:general}.  Item (2) of Theorem 
\ref{theorem:rigidity} follows from item (3) together with the theorem of Suslin-Abe that 
for these rings $R$, we have $G(R)=E(R)$.
\endproof

\subsection{Reflection groups}
 
Let $\Gamma_n$ be a discrete group generated by reflections in the faces 
of a compact $n$-simplex $P$ in the space $Y=\R^n$ or 
$Y=\H^n$.  As $\Gamma_n$ is discrete, the subgroup $H_v$ of $\Gamma_n$
stabilizing any one of the $n+1$ vertices $v_i$ of $P$ is a finite group.  Note that 
the groups $\{H_v:v\in P\}$ generate all of $\Gamma_n$.  Any proper 
subcollection of stabilizers of the 
$v_i$ generates a group which stabilizes some face of $P$, and is therefore finite. 

Now suppose that $\Gamma_n$ acts on an $m$-dimensional space $X$ with 
$m\leq n-1$, satisfying the axioms in the definition of strong property $\FA_m$.  Then by the previous paragraph,  and by the Bruhat-Tits Fixed Point Theorem (Lemma \ref{lemma:finite} above), each stabilizer of a 
vertex $v_i$ in $P$ has nonempty fixed set, as does the intersection of the stabilizers of the 
vertices in any codimension one face of $P$.  It follows that the nerve of the cover of fixed sets 
is an $n$-simplex.  Since $\dim(X)=m<n$, it follows just as in the proof of 
Theorem \ref{theorem:general}  that the intersection of these fixed sets is nonempty, and we are done.
\endproof

\bigskip

\noindent
Benson Farb:\\
Department of Mathematics, Univeristy of Chicago\\
5734 University Ave.\\
Chicago, Il 60637\\
E-mail: farb@math.uchicago.edu

\end{document}